\documentclass[10pt,reqno]{amsart}
\usepackage{amsmath}
\usepackage{amssymb}
\usepackage{amsthm}

\newlength{\defbaselineskip}
\newcommand{\setlinespacing}[1]%
           {\setlength{\baselineskip}{#1 \defbaselineskip}}

\numberwithin{equation}{section}

\newtheorem{thm}{Theorem}[section]

\newtheorem{prop}[thm]{Proposition}

\theoremstyle{definition}

\theoremstyle{remark}
\newtheorem{rem}[thm]{Remark}
\numberwithin{equation}{section}

\begin{document}

\title[The Schr\"odinger smoothing effect]
{A note on the Schr\"odinger smoothing effect}

\author{Ihyeok Seo}

\thanks{I. Seo was supported by the NRF grant funded by the Korea government(MSIP) (No. 2017R1C1B5017496).}

\subjclass[2010]{Primary: 35B45; Secondary: 42B37}
\keywords{Smoothing estimates, Fourier restriction.}

\address{Department of Mathematics, Sungkyunkwan University, Suwon 16419, Republic of Korea}
\email{ihseo@skku.edu}

\begin{abstract}
The Kato-Yajima smoothing estimate is a smoothing weighted $L^2$ estimate with a singular power weight
for the Schr\"odinger propagator. The weight has been generalized relatively recently to Morrey-Campanato weights.
In this paper we make this generalization more sharp in terms of the so-called Kerman-Sawyer weights.
Our result is based on a more sharpened Fourier restriction estimate in a weighted $L^2$ space.
Obtained results are also extended to the fractional Schr\"odinger propagator.
\end{abstract}

\maketitle

\section{Introduction}

Consider the following free Schr\"odinger equation which describes the wave behavior of a quantum particle
that is moving in the absence of external forces:
\begin{equation*}
\begin{cases}
i\partial_tu+\Delta u=0,\\
u(x,0)=f(x).
\end{cases}
\end{equation*}
Applying the Fourier transform to the equation, the solution is then given by
\begin{equation*}
e^{it\Delta}f(x) :=(2\pi)^{-n}\int_{\mathbb{R}^n}e^{ix\cdot\xi}e^{-it|\xi|^2} \widehat{f}(\xi)d\xi.
\end{equation*}

The physical interpretation for the solution is that
$|u(x,t)|^2$ is the probability density for finding the particle at place $x\in\mathbb{R}^n$ and time $t\in\mathbb{R}$.
Indeed,
$\|e^{it\Delta}f\|_{L_x^2}=\|f\|_{L^2}$ for any fixed $t$.
But interestingly, if we take $L^2$-norm in $t$,
the extra gain of regularity of order $1/2$ in $x$ can be observed in the case $n=1$.
That is,
for any fixed $x\in\mathbb{R}$,
$$\||\nabla|^{1/2}e^{it\Delta}f\|_{L_t^2}=C\|f\|_{L^2}.$$
These facts follow easily from Plancherel's theorem.

This kind of smoothing effect for higher dimensions
is known to be valid in a weighted $L^2$ space with a singular power weight in the spatial variable.
In fact,
\begin{equation}\label{KYest}
\||\nabla|^{s}e^{it\Delta}f\|_{L_{x,t}^2(|x|^{-2(1-s)})}\leq C\|f\|_{L^2}
\end{equation}
if and only if $-\frac{n-2}2<s<\frac12$ and $n\geq2$.
These so called Kato-Yajima smoothing estimates were first discovered by Kato and Yajima \cite{KY}
for $0\leq s<\frac12$ when $n\geq3$, and $0<s<\frac12$ for $n=2$.
Ben-Artzi and Klainerman \cite{BK} gave an alternate proof of this result.
Since then, Watanabe \cite{W} showed that \eqref{KYest} fails for the case $s=\frac12$,
and the full range was obtained by Sugimoto \cite{Su} although it was later shown by Vilela \cite{V} that the range is indeed optimal.
In addition, it was shown in \cite{RS} that the critical case $s=1/2$ can be attained under some structure condition on the operator $|x|^{-1/2}|\nabla|^{1/2}$.

Now we discuss a generalization of \eqref{KYest} to Morrey-Campanato weights $w\in\mathcal{L}^{\alpha,p}$
defined by
$$\|w\|_{\mathcal{L}^{\alpha,p}}
:=\sup_{x\in\mathbb{R}^{n},r>0}r^\alpha
\bigg( \frac{1}{r^{n}} \int_{B(x,r)} w(y)^p dy\bigg)^{1/p}<\infty$$
for $\alpha>0$ and $1\leq p\leq n/\alpha$.
Here, $B(x,r)$ denotes a ball in $\mathbb{R}^{n}$ centered at $x$ with radius $r$.
Notice that $\mathcal{L}^{\alpha,p} =L^p$ when $p=n/\alpha$,
and even $L^{n/\alpha,\infty}\subset\mathcal{L}^{\alpha,p}$ for $p<n/\alpha$.
So, the class $\mathcal{L}^{\alpha,p}$ is a natural extension of $L^p$ and contains much singular functions
like $|x|^{-\alpha}$ which does not belong to any $L^p$ class.
By observing that $|x|^{-2(1-s)}\in\mathcal{L}^{2(1-s),p}$ for $p<\frac n{2(1-s)}$
together with \eqref{KYest},
one can naturally expect a more general estimate,
\begin{equation}\label{general}
\||\nabla|^se^{it\Delta}f\|_{L_{x,t}^2(w(x))}\leq C\|w\|_{\mathcal{L}^{2(1-s),p}}^{1/2}\|f\|_{L^2},
\end{equation}
for some $-\frac{n-2}{2}< s<\frac12$ and $1\leq p< \frac n{2(1-s)}$.
Indeed, Barcel\'{o} et al \cite{BBCRV} obtained this generalized smoothing estimate \eqref{general}\footnote{Estimate \eqref{general}
was also shown in \cite{BBCRV} for the border line $p=\frac n{2(1-s)}$, but in this case $\mathcal{L}^{2(1-s),p}=L^p$.
In view of generalization of \eqref{KYest} with a singular power weight, we therefore exclude it in \eqref{general}.} for
\begin{equation}\label{range}
-\frac{n-2}{2}< s<\frac1{n+1}\quad \text{and}\quad \frac{n-1}{2(1-2s)}<p< \frac n{2(1-s)}.
\end{equation}
(See also \cite{RV} for the case $s=0$.)
In the same paper they also obtained some necessary conditions which shows that \eqref{range} is sharp particularly when $s>0$.
In other words, the smoothing effect with a gain of regularity of order $s\geq\frac1{n+1}$ is not possible.

In this paper we make the estimate \eqref{general} more sharp.
To present results, we first introduce the Kerman-Sawyer class $\mathcal{KS}_{\alpha}$ which is defined for $0<\alpha<n$\, if
$$\|w\|_{\mathcal{KS}_{\alpha}} := \sup_Q \left( \int_Q w(x) dx\right)^{-1} \int_Q \int_Q \frac{w(x)w(y)}{|x-y|^{n-\alpha}}  dxdy < \infty,$$
where the sup is taken over all dyadic cubes $Q$ in $\mathbb{R}^n$.
This class of weights gives the characterization of weighted $L^2$ estimates for the fractional integrals.
Notice that $\mathcal{KS}_{\alpha}\subset\mathcal{L}^{\alpha,1}$, but in general,
\begin{equation}\label{wider}
\mathcal{L}^{\alpha, p} \subset \mathcal{KS}_{\alpha}\quad\text{if}\quad p\neq1.
\end{equation}
For more details we refer the reader to \cite{KS} (see Theorem 2.3 there).
See also \cite{BBRV}, Subsection 2.2.

Now we state our main result which makes the estimate \eqref{general} more sharp
in terms of the size $\|w^\beta\|^{1/\beta}_{\mathcal{KS}_{\alpha}}$ instead of $\|w\|_{\mathcal{L}^{2(1-s),p}}$
(see Remark \ref{rem} below).

\begin{thm}\label{thm}
Let $n\geq2$. If\, $1\leq\beta<\frac{n+1}2$ and $\alpha>\beta+\frac{n-1}2$, then
\begin{equation}\label{improve}
\||\nabla|^se^{it\Delta}f\|_{L_{x,t}^2(w(x))}\leq C\|w^\beta\|_{\mathcal{KS}_\alpha}^{\frac1{2\beta}}\|f\|_{L^2}
\end{equation}
with\, $s=1-\alpha/(2\beta)$.
\end{thm}

\begin{rem}
Inserting $\beta=\frac{\alpha}{2(1-s)}$ into the conditions $1\leq\beta<\frac{n+1}2$ and $\alpha>\beta+\frac{n-1}2$, one can see that
$-\frac{\alpha-2}{2}\leq s<\frac{1-(n-\alpha)}{n+1-2(n-\alpha)}$ in the theorem.
Since $0<\alpha<n$, when $\alpha\rightarrow n$,
our theorem therefore recovers the optimal range of $s$ in \eqref{range}.
\end{rem}

\begin{rem}\label{rem}
From \eqref{wider}, if $1<q<n/\alpha$,
$\|w^\beta\|^{1/\beta}_{\mathcal{KS}_\alpha}\leq\|w^\beta\|^{1/\beta}_{\mathcal{L}^{\alpha,q}}=\|w\|_{\mathcal{L}^{\alpha/\beta,q\beta}}.$
Using this fact and setting $\alpha/\beta=2(1-s)$ and $q\beta=p$,
one can easily check that
the conditions in the theorem together with $1<q<n/\alpha$ imply
$$-\frac{n-2}{2}< s<\frac1{n+1}\quad\text{and}\quad\frac{n-1}{2(1-2s)}<\beta<p<\frac{n}{2(1-s)}.$$
Hence \eqref{improve} implies \eqref{general} under the same conditions on
$s$ and $p$ as in \eqref{range}.
\end{rem}

The quantity $\|f^\beta\|^{1/\beta}_{\mathcal{KS}_\alpha}$ was already appeared in \cite{S} and \cite{LS} concerning
the unique continuation for the Schr\"odinger equation and eigenvalue bounds for the Schr\"odinger operator, respectively.
We also refer the reader to \cite{BBCRV2,KoS} for some weighted $L^2$ estimates in which a time-dependent weight $w(x,t)$ is involved.

The study of smoothing estimates is particularly important in applications to nonlinear dispersive equations
especially with derivatives in the potential or in the nonlinearity.
Since 1980's spectral methods based on resolvent estimates were developed to understand these smoothing effects for a variety of dispersive equations,
starting from Kato's work \cite{K},
and since 1990's harmonic analysis methods have been also used.
Recently, a completely different new approach was also
developed by Ruzhansky and Sugimoto \cite{RS2} using canonical transforms and methods of comparison.
Here we follow the harmonic analysis approach which bypasses resolvent estimates to obtain more directly smoothing estimates.
The key ingredient in this approach is the availability of a weighted $L^2$ restriction estimate for the Fourier transform.
In our case we need to obtain the following restriction estimate:

\begin{prop}\label{prop}
Let $n\geq2$. If $1\leq\beta<\frac{n+1}2$ and $\alpha>\beta+\frac{n-1}2$, then
\begin{equation}\label{restri}
\|\widehat{fd\sigma}\|_{L^2(w)}\leq C\|w^\beta\|_{\mathcal{KS}_\alpha}^{\frac1{2\beta}}\|f\|_{L^2(\mathbb{S}^{n-1})}.
\end{equation}
\end{prop}

The Stein-Tomas restriction theorem (see \cite{MS}), in its dual form, states that
\begin{equation*}
\|\widehat{fd\sigma}\|_{L^r}\leq C\|f\|_{L^2(\mathbb{S}^{n-1})}
\end{equation*}
for $r\geq2(n+1)/(n-1)$.
Applying H\"older's inequality, this estimate now implies
\begin{equation}\label{restri2}
\|\widehat{fd\sigma}\|_{L^2(w)}\leq C\|w\|_{L^p}^{1/2}\|f\|_{L^2(\mathbb{S}^{n-1})}
\end{equation}
whenever $1\leq p\leq\frac{n+1}2$.
In \cite{RV2}, Ruiz and Vega extended \eqref{restri2} by replacing the $L^p$ norm of $w$ by the Morrey-Campanato norm $L^{\alpha,p}$
with $\frac{2n}{n+1}<\alpha\leq n$ and $\frac{n-1}{2(\alpha-1)}<p<\frac{n}{\alpha}$
(see also \cite{CS,CR,RV} for $\alpha=2$).
Estimate \eqref{restri} makes this extension more sharp as in Remark \ref{rem}.
See also \cite{S} for the case where $\alpha=n-1$ and $\beta=\frac{n-1}{2}$ in \eqref{restri},
which just implies
\begin{equation*}
\|e^{it\Delta}f\|_{L_{x,t}^2(w(x))}\leq C\|w^{\frac{n-1}{2}}\|_{\mathcal{KS}_{n-1}}^{\frac1{n-1}}\|f\|_{L^2}.
\end{equation*}

We discuss further the smoothing estimate for the fractional Schr\"odinger equation
which has recently attracted interest from mathematical physics.
By generalizing the Feynman path integral to the L\'{e}vy one,
Laskin \cite{L} introduced the fractional quantum mechanics in which
the following fractional Schr\"odinger equation plays a central role
and it is conjectured that physical realizations may be limited to $1<\gamma<2$:
\begin{equation*}
\begin{cases}
i\partial_tu-(-\Delta)^{\gamma/2} u=0,\\
u(x,0)=f(x),
\end{cases}
\end{equation*}
where $(-\Delta)^{\gamma/2}$ is defined by means of the Fourier transform:
$$\widehat{(-\Delta)^{\gamma/2}f}(\xi)=|\xi|^\gamma\widehat{f}(\xi).$$

Applying the Fourier transform to the equation as before, the solution is given by
\begin{equation}\label{erty}
e^{-it(-\Delta)^{\gamma/2}}f(x) :=(2\pi)^{-n}\int_{\mathbb{R}^n}e^{ix\cdot\xi}e^{-it|\xi|^\gamma} \widehat{f}(\xi)d\xi.
\end{equation}
In \cite{W}, Watanabe extended the Kato-Yajima smoothing estimate \eqref{KYest} to the fractional case \eqref{erty}, as follows:
\begin{equation}\label{wata}
\||\nabla|^{s}e^{-it(-\Delta)^{\gamma/2}}f\|_{L_{x,t}^2(|x|^{-(\gamma-2s)})}\leq C\|f\|_{L^2}
\end{equation}
for $\gamma>1$, $s\geq0$, $-\frac{n-\gamma}2<s<\frac{\gamma-1}2$, and $n\geq2$.
In the following theorem, which recovers Theorem \ref{thm} when $\gamma=2$,
we generalize \eqref{wata} in the spirit of the Schr\"odinger case $\gamma=2$.

\begin{thm}\label{thm0}
Let $n\geq2$ and $\gamma>1$. If\, $1\leq\beta<\frac{n+1}2$ and $\alpha>\beta+\frac{n-1}2$, then
\begin{equation}\label{frac}
\||\nabla|^se^{-it(-\Delta)^{\gamma/2}}f\|_{L_{x,t}^2(w(x))}\leq C\|w^\beta\|_{\mathcal{KS}_\alpha}^{\frac1{2\beta}}\|f\|_{L^2}
\end{equation}
with
\begin{equation}\label{condition}
s=1-\frac12(\frac\alpha\beta+2-\gamma).
\end{equation}
\end{thm}

\begin{rem}
Particularly when $1<\gamma\leq2$, from applying the conditions $\alpha>\beta+\frac{n-1}2$ and $1\leq\beta<\frac{n+1}2$ in turn to \eqref{condition},
it follows that
\begin{align*}
s&<1-\frac12(1+\frac{n-1}{2\beta}+2-\gamma)\\
&<1-\frac12(1+\frac{n-1}{n+1}+2-\gamma)=\frac{1}{n+1}-\frac12(2-\gamma).
\end{align*}
Notice here that this range when $\gamma=2$ is reduced to the Schr\"odinger case,
and that it implies the restriction $\gamma>\frac{6n+2}{3(n+1)}$ when $s>0$.
Hence \eqref{frac} gives a smoothing effect for the fractional Schr\"odinger propagator $e^{-it(-\Delta)^{\gamma/2}}$
of order $\gamma>\frac{2n}{n+1}$.
\end{rem}

Throughout this paper, the letter $C$ stands for a positive constant which may be different
at each occurrence.
We also denote $A\lesssim B$ and $A\sim B$ to mean $A\leq CB$ and $CB\leq A\leq CB$, respectively,
with unspecified constants $C>0$.

\section{Proofs}

\subsection{Proof of Theorem \ref{thm0}}
Using polar coordinates ($\xi\rightarrow\rho\sigma$) and a change of variables ($\rho^\gamma\rightarrow r$),
we first write
\begin{align*}
|\nabla|^se^{-it(-\Delta)^{\gamma/2}}f(x)&=\int_{\mathbb{R}^n} e^{ix\cdot\xi}|\xi|^se^{-it|\xi|^\gamma}\widehat{f}(\xi)d\xi\\
&=\int_0^\infty\int_{\mathbb{S}^{n-1}}e^{ix\cdot\rho\sigma}\rho^se^{-it\rho^\gamma}\widehat{f}(\rho\sigma)\rho^{n-1}d\sigma d\rho\\
&\sim\int_0^\infty\int_{\mathbb{S}_{r^{1/\gamma}}^{n-1}}e^{ix\cdot r^{1/\gamma}\sigma}r^{s/\gamma}e^{-itr}
\widehat{f}(r^{1/\gamma}\sigma)d\sigma_{r^{1/\gamma}}r^{-(1-1/\gamma)}dr\\
&=\int_{\mathbb{R}} e^{-itr}\chi_{(0,\infty)}(r)r^{\frac{s+1}{\gamma}-1} \widehat{\widehat{f}d\sigma_{r^{1/\gamma}}}(-x)dr.
\end{align*}
Then we apply Plancherel's theorem in the $t$-variable to get
\begin{align*}
\||\nabla|^se^{-it(-\Delta)^{\gamma/2}}f\|_{L_{x,t}^2(w(x))}
&\sim\|\chi_{(0,\infty)}(r)r^{\frac{s+1}{\gamma}-1}\widehat{\widehat{f}d\sigma_{r^{1/\gamma}}}(-x)\|_{L_{x,r}^2(w(x))}\\
&=\|\chi_{(0,\infty)}(r)r^{\frac{s+1}{\gamma}-1}\|\widehat{\widehat{f}d\sigma_{r^{1/\gamma}}}(-x)\|_{L_x^2(w(x))}\|_{L_r^2}.
\end{align*}

Now we use the following re-scaled estimate of \eqref{restri},
$$\|\widehat{fd\sigma_{r^{1/\gamma}}}\|_{L^2(w)}\leq Cr^{\frac1{2\gamma}(\frac\alpha\beta-1)}\|w^\beta\|_{\mathcal{KS}_\alpha}^{\frac1{2\beta}}\|f\|_{L^2(\mathbb{S}_{r^{1/\gamma}}^{n-1})}$$
for $1\leq\beta<\frac{n+1}2$ and $\alpha>\beta+\frac{n-1}2$.
Then,
\begin{align*}
\||\nabla|^se^{-it(-\Delta)^{\gamma/2}}f\|_{L_{x,t}^2(w(x))}
&\lesssim\|w^\beta\|_{\mathcal{KS}_\alpha}^{\frac1{2\beta}}\|\chi_{(0,\infty)}(r)r^{\frac{s+1}{\gamma}-1} r^{\frac1{2\gamma}(\frac\alpha\beta-1)}\|\widehat{f}\|_{L^2(\mathbb{S}_{r^{1/\gamma}}^{n-1})}\|_{L_r^2}\\
&=\|w^\beta\|_{\mathcal{KS}_\alpha}^{\frac1{2\beta}}\bigg(\int_0^\infty r^{2(\frac{s+1}{\gamma}-1)+\frac1{\gamma}(\frac{\alpha}{\beta}-1)}
\int_{\mathbb{S}_{r^{1/\gamma}}}|\widehat{f}|^2d\sigma_{r^{1/\gamma}}dr\bigg)^{\frac12}\\
&\sim\|w^\beta\|_{\mathcal{KS}_\alpha}^{\frac1{2\beta}}\|\widehat{f}\|_{L^2}
\end{align*}
if $2(\frac{s+1}{\gamma}-1)+\frac1{\gamma}(\frac{\alpha}{\beta}-1)=\frac1\gamma-1$
which is equivalent to the condition \eqref{condition}.
Finally, Plancherel's theorem gives the desired estimate \eqref{frac}
for $1\leq\beta<\frac{n+1}2$ and $\alpha>\beta+\frac{n-1}2$
with $s=1-\frac12(\frac\alpha\beta+2-\gamma)$.

\subsection{Proof of Proposition \ref{prop}}

By the standard $TT^\ast$ argument it is enough to show that
\begin{equation}\label{dl}
\|\widehat{d\sigma}\ast f\|_{L^2(w)}\leq C\|w^\beta\|_{\mathcal{KS}_\alpha}^{\frac1{\beta}}\|f\|_{L^2(w^{-1})}.
\end{equation}
Indeed, using \eqref{dl}, we see that
\begin{align*}
\int_{\mathbb{S}^{n-1}}|\widehat{f}|^2d\sigma=\int_{\mathbb{S}^{n-1}}\widehat{f}\overline{\widehat{f}}d\sigma
&=\int_{\mathbb{R}^n}f(f\ast\widehat{d\sigma})dx\\
&\leq\|f\|_{L^2(w^{-1})}\|f\ast\widehat{d\sigma}\|_{L^2(w)}\\
&\leq C\|w^\beta\|_{\mathcal{KS}_\alpha}^{\frac1{\beta}}\|f\|_{L^2(w^{-1})}^2.
\end{align*}
Namely, we get
$$\|\widehat{f}\|_{L^2(\mathbb{S}^{n-1})}\leq \|w^\beta\|_{\mathcal{KS}_\alpha}^{\frac1{2\beta}}\|f\|_{L^2(w^{-1})}$$
which is equivalent to
$$\|\widehat{fd\sigma}\|_{L^2(w)}\leq \|w^\beta\|_{\mathcal{KS}_\alpha}^{\frac1{2\beta}}\|f\|_{L^2(\mathbb{S}^{n-1})}$$
by duality.
(The first estimate \eqref{dl} can be also deduced immediately from the last two equivalent estimates,
and so all these estimates are equivalent each other.)

To show \eqref{dl}, we first decompose $\widehat{d\sigma}(x)$ dyadically as
\begin{align*}
\widehat{d\sigma}(x)&=\sum_{j=0}^\infty \widehat{d\sigma}(x)\psi_j(|x|)\\
&:=\sum_{j=0}^\infty K_j(x),
\end{align*}
where supp\,$\psi_0\subset[-1,1]$ and supp\,$\psi_j\subset[2^{j-1},2^{j+1}]$ for $j\geq1$.
Then we get \eqref{dl} by interpolating the two estimates
\begin{equation}\label{dl2}
\|K_j\ast f\|_{L^2}\leq C2^j\|f\|_{L^2}
\end{equation}
and
\begin{equation}\label{dl3}
\|K_j\ast f\|_{L^2(w^\beta)}\leq C2^{j(n-\alpha-(n-1)/2)}\|w^\beta\|_{\mathcal{KS}_\alpha}\|f\|_{L^2(w^{-\beta})}.
\end{equation}
Indeed, by the interpolation, it follows that
\begin{equation*}
\|K_j\ast f\|_{(L^2,L^2(w^\beta))_{1/\beta,2}}\leq C2^{j(1+\frac1\beta(\frac{n-1}2-\alpha))}\|w^\beta\|_{\mathcal{KS}_\alpha}^{1/\beta}\|f\|_{(L^2,L^2(w^{-\beta}))_{1/\beta,2}}.
\end{equation*}
Now we use the following real interpolation space identity (see Theorem 5.4.1 in \cite{BL}),
$$( L^{2} (w_0), L^{2} (w_1) )_{\theta,2} = L^2(w),\quad w= w_0^{1-\theta} w_1^{\theta},$$
with $0<\theta=1/\beta\leq1$, $w_0=1$, and $w_1=w^\beta \text{ or } w^{-\beta}$,
and then summing on $j$, we get \eqref{dl}
if $1+\frac1\beta(\frac{n-1}2-\alpha)<0$ which is equivalent to the condition $\alpha>\beta+\frac{n-1}{2}$ in the theorem.

It remains to show \eqref{dl2} and \eqref{dl3}.
The first estimate \eqref{dl2} was already well-known in the proof of the Stein-Tomas restriction theorem (see \cite{MS}).
On the other hand, the argument in \cite{RV2} can be adapted for the second estimate \eqref{dl3}.
Indeed, we first make in $\mathbb{R}^n$ a grid with cubes $Q_\nu$ with side length $2^j$,
and set $f_\nu=f\chi_{Q_\nu}$. Then,
\begin{align}\label{11}
\nonumber\int_{\mathbb{R}^n}|K_j\ast f|^2w^\beta dx&=\int_{\mathbb{R}^n}|K_j\ast \sum_\nu f_\nu|^2w^\beta dx\\
\nonumber&\lesssim\sum_\nu\int_{\mathbb{R}^n}|K_j\ast f_\nu|^2w^\beta dx\\
&\lesssim\bigg(\sup_\nu\int_{Q_\nu^\ast}w^\beta dx\bigg)\sum_\nu\|K_j\ast f_\nu\|_{L^\infty(Q_\nu^\ast)}^2,
\end{align}
where $Q_\nu^\ast$ is a cube with the same center as $Q_\nu$ and side ten times bigger than the sides of $Q_\nu$.
Now we use the following well-known decay estimate (see \cite{MS}) for the Fourier transform of the surface measure $d\sigma$,
$$|\widehat{d\sigma}(x)|\leq C(1+|x|)^{-(n-1)/2},$$
to get
\begin{align}\label{22}
\nonumber\sum_\nu\|K_j\ast f_\nu\|_{L^\infty(Q_\nu^\ast)}^2
&\leq\sum_\nu\|K_j\|_{L^\infty(Q_\nu^\ast)}^2\|f_\nu\|_{L^1}^2\\
\nonumber&\lesssim2^{-j(n-1)}\sum_\nu\|f_\nu\|_{L^1}^2\\
\nonumber&\lesssim2^{-j(n-1)}\sum_\nu\|w^{\beta/2}\|_{L^2(Q_\nu)}^2\|w^{-\beta/2}f\|_{L^2(Q_\nu)}^2\\
\nonumber&\lesssim2^{-j(n-1)}\bigg(\sup_\nu\int_{Q_\nu}w^\beta dx\bigg)\sum_\nu\|w^{-\beta/2}f\|_{L^2(Q_\nu)}^2\\
&\lesssim2^{-j(n-1)}\bigg(\sup_\nu\int_{Q_\nu}w^\beta dx\bigg)\|f\|_{L^2(w^{-\beta})}^2.
\end{align}
Combining \eqref{11} and \eqref{22}, we get the desired estimate
\begin{align*}
\|K_j\ast f\|_{L^2(w^\beta)}&\leq C2^{-j(n-1)/2}\bigg(\sup_\nu\int_{Q_\nu^\ast}w^\beta dx\bigg)\|f\|_{L^2(w^{-\beta})}\\
&\lesssim2^{-j(n-1)/2}2^{j(n-\alpha)}\|w^\beta\|_{\mathcal{KS}_{\alpha}}\|f\|_{L^2(w^{-\beta})}.
\end{align*}

\

\noindent\textbf{Acknowledgment.}
The author thanks the anonymous referees for their careful reading of his manuscript
and their helpful comments.

\end{document}